# Algorithms for Learning and Teaching Sets of Vertices in Graphs

Patricia A. Evans and Michael R. Fellows
Department of Computer Science
University of Victoria
Victoria, B.C. V8W 3P6, Canada

Lane H. Clark and Walter D. Wallis
Department of Mathematics
Southern Illinois University
Carbondale, Illinois 62901-4408, U.S.A.

## Abstract

The learning complexity of special sets of vertices in graphs is studied in the model(s) of *exact learning by (extended) equivalence and membership queries*. Polynomial-time learning algorithms are described for vertex covers, independent sets and dominating sets. The complexity of learning vertex sets of fixed size is also investigated, and it is shown that the $k$-element vertex covers in a graph can be learned in a number of rounds of interaction that it is *independent of the size of the graph*. Apart from the elegance of these algorithmic problems, the chief motivation is the surprising recently established connection between the important unsolved problem of the learning complexity of CNF (or DNF) formulas and the learning complexity of dominating sets. The complexity of *teaching* sets of vertices in graphs is also considered.

## 1 Introduction and Preliminaries

The complexity of *learning* combinatorial objects of various kinds, such as finite automata [Ang87] and Boolean formulas [Val84] has recently become an active area of research in computational complexity. There is now an annual conference series in computational learning theory (COLT) and important connections between research in this area and other subjects, such as cryptography and structural complexity theory. In the context of this expanding program of research it would seem to be natural to study the learning complexity of simple structures in graphs. One of the contributions of this paper is to present some fundamental results on learning algorithms in this hitherto unexplored setting.

*There is, however, a far more important motivation for considering the learning complexity of graph structures.* Since its inception, one of the most important unsolved problems in computational learning theory has been whether it is possible to learn the truth assignments to a CNF (or DNF) formula in polynomial time [Val84]. (See [Ang92] for a survey of what is presently known.) In [DEF93] the following connections to the learning complexity of graph structures are established.

*Theorem 1.1* [DEF93]. In the model of exact learning by extended equivalence queries, the truth assignments of a CNF formula can be learned in polynomial time if and only the dominating sets of a graph can be learned in polynomial time.

*Theorem 1.2* [DEF93]. In the model of exact learning by extended equivalence and membership queries, the truth assignments of a CNF formula can be learned in polynomial time if and only if the there is a polynomial $q$ and a $k$-uniform learning algorithm $L$ that learns the $k$-element dominating sets in a graph of order $n$ in time $q(n)$.

It is probably fair to say that most researchers in computational learning theory favor the conjecture that CNF is not polynomial-time learnable. Assuming that this is the case, and in view of the positive results which we obtain concerning vertex covers and independent sets, it would seem that the boundary between what can and what cannot be learned in polynomial time can be investigated in an interesting way in the graph-theoretic setting. The boundary seems to be somewhere "between" independent sets and dominating sets.

As suggested by the above theorems we study the complexity of learning, e.g., (1) *all* the vertex covers in a graph, and (2) just the vertex covers of a fixed size $k$.

The model(s) of learning we employ are essentially those of Angluin [Ang87], [Ang92].

In the model of *exact learning by equivalence (and membership) queries* we consider that there are two players, the Learner and the Teacher. We assume that the Teacher possesses privately a graph $G = (V, E)$ with vertex set $V = \{1, \ldots, n\}$, and that the Learner initially knows only the vertex set, that is, essentially just the order $n$ of the graph. The goal of the interaction which follows is for the Learner to produce a graph on this vertex set which is *equivalent* to the one possessed by the Teacher with respect to the kinds of vertex sets which are being taught.

*Example. Learning Dominating Sets.*
The goal of the Learner in this case is to produce a graph $H$ that is equivalent to the graph $G$ of the Teacher in the sense that for any set $V' \subseteq V = \{1, \ldots, n\}$, $V'$ is a dominating set in $H$ if and only if $V'$ is a dominating set in $G$. The Learner will then possess a representation of the property of a set of vertices of being a dominating set in $G$.

*Definition.* Let $\Phi(G, V')$ be a predicate representing a property of vertex sets $V'$ of finite graphs $G$. Write $S(G) = \{V' : \Phi(G, V')\}$ and term $S$ the set of *solutions* for $G$. Graphs $G$ and $H$ are $\Phi$-*equivalent* if (1) $V(G) = V(H)$, and (2) $S(G) = S(H)$. We will term such a predicate $\Phi$ a *vertex set concept*.

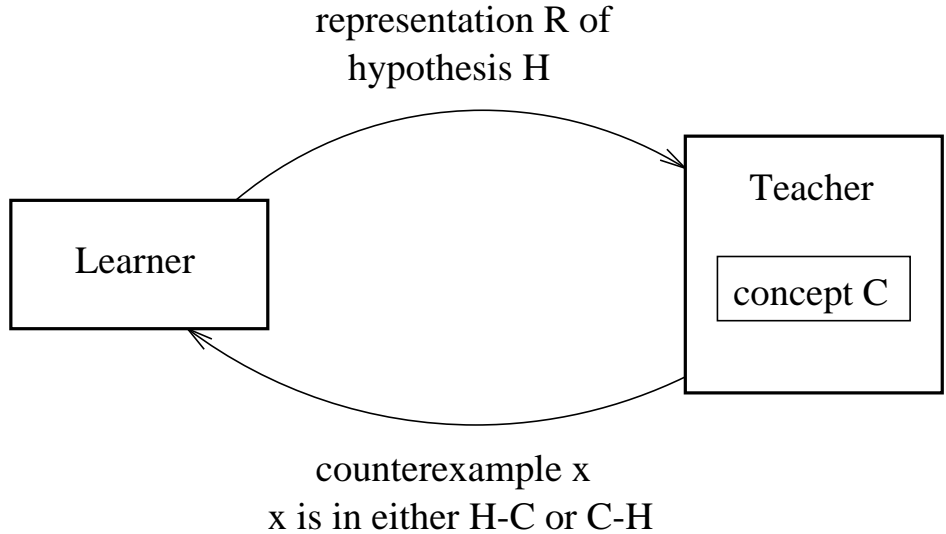

Figure 1: Interaction between Learner and Teacher in an equivalence query.

There are two kinds of interaction between the Learner and the Teacher that we will consider:

(1) A *membership query* to the Teacher consists in the presentation to the Teacher of a set of vertices $V' \subseteq \{1, \ldots, n\}$. The response of the Teacher is either "yes" or "no" depending on whether $\Phi(G, V')$.

(2) An *equivalence query* to the Teacher consists in the presentation to the Teacher of a graph $H$ on the vertex set $\{1, \ldots, n\}$. (We may refer to this graph as the *hypothesis* of the Learner.) The response of the Teacher is the message, "finished," in the case that $H$ is equivalent to the graph $G$ being taught, or a message giving a counterexample in the case where $H$ and $G$ are not equivalent. Figure 1 illustrates this interaction between the Learner and Teacher. A *counterexample* is either:

(a) A set of vertices $V'$ such that $\Phi(G, V')$ holds, but not $\Phi(H, V')$, which we term a *positive counterexample*, or

(b) A set of vertices $V'$ such that $\Phi(H, V')$ holds, but not $\Phi(G, V')$, which we term a *negative counterexample*.

A *learning algorithm* $L$ for a predicate $\Phi$ is a strategy by which the Learner can always produce a graph $H$ that is $\Phi$-equivalent to the graph $G$ possessed by the Teacher, assuming only that the Teacher is *competent*, that is, responds properly to the queries made by the

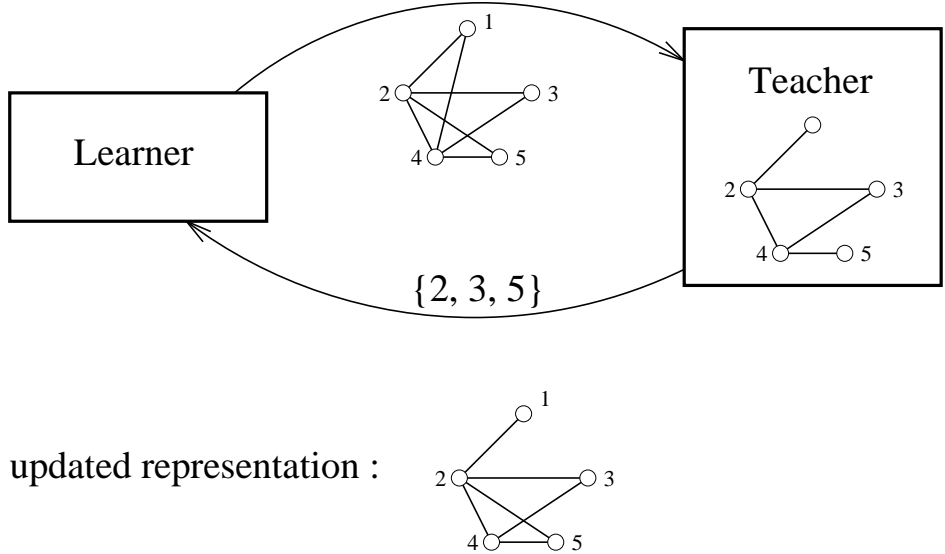

Figure 2: Update cycle for learning Vertex Covers by Equivalence Queries.

Learner. The learning algorithm $L$ is said to run in $f(n)$ *rounds* if for any graph $G$ of order $n$ that is being taught by a competent but not necessarily optimal Teacher, the algorithm finishes after making at most $f(n)$ queries to the Teacher.

In the model of *exact learning by* extended *equivalence (and membership) queries* the framework is the same, except that instead of offering a graph $H$ on $\{1,...,n\}$ as a hypothesis, the Learner may present a circuit $C$ with $n$ inputs and a single output (a decision circuit). The circuit $C$ is *equivalent* to $G$ if it has output 1 for precisely the 0-1 vectors which correspond to the vertex sets being taught. (For a set $U \subset \{1,...,n\}$, the corresponding vector has $j^{th}$ component 1 if and only if $j \in U$.)

## 2  Learning Vertex Sets

We consider in this section learning algorithms for the following vertex set concepts.
(1) $V'$ is a *vertex cover* in $G$: for every edge $uv$ of $G$, either $u \in V'$ or $v \in V'$.
(2) $V'$ is an *independent set* in $G$: for every pair of vertices $u, v \in V'$, $uv$ is not an edge of $G$.
(3) $V'$ is a *dominating set* in $G$: for every vertex $u$, either $u \in V'$ or $u$ is adjacent to a vertex $v$ with $v \in V'$.

We will say that $V'$ is a $k$-vertex cover in $G$ (and similarly for other kinds of vertex sets) if $V'$ is a vertex cover in $G$ of cardinality $k$.

*Theorem 2.1* Vertex Covers can be learned by equivalence queries in $O(n^2)$ rounds.

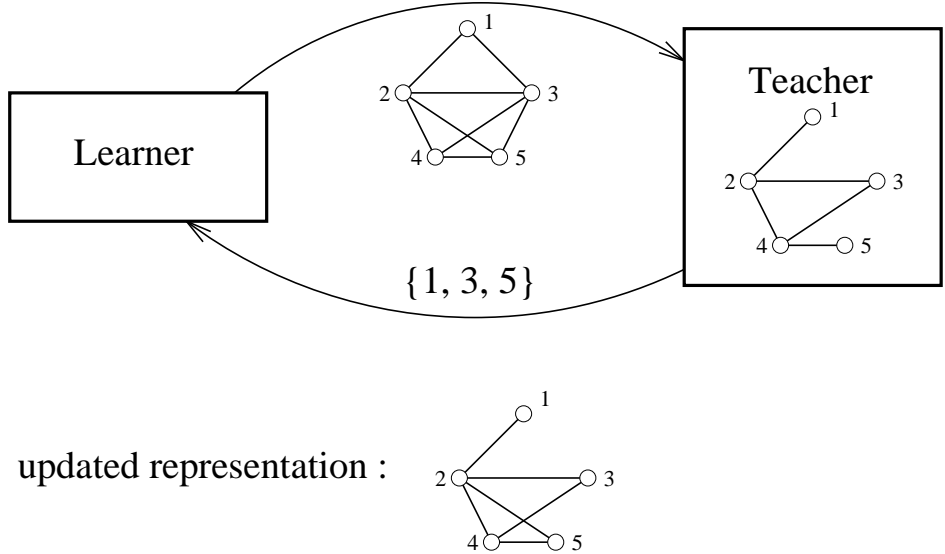

Figure 3: Update cycle for learning Independent Sets by Equivalence Queries.

*Proof.* The Learner begins with an equivalence query consisting of the complete graph. If $G$ is not complete, then the Teacher must respond with a positive counterexample $V'$ such that $|V - V'| \geq 2$. The Learner can deduce from this counterexample that for any two vertices $u, v$ of $V - V'$, $uv$ is not an edge of $G$. The Learner continues to offer hypotheses according to the following policy: the hypothesis graph $H$ contains all possible edges except those for which the Learner has deduced in the above fashion that the edge is not present in $G$. The edge set of the hypothesis graph $H$ is therefore always a superset of the edge set of $G$. From this it follows that every response of the Teacher to an equivalence query (the only kind made by this algorithm) must consist of a positive counterexample. An example of a query and counterexample round is shown in Figure 2. We may also observe that any positive counterexample must provide an opportunity for the Learner to deduce the non-presence of at least one new edge in $G$. Thus the algorithm will finish correctly in at most $\binom{n}{2}$ rounds. □

*Theorem 2.2* Independent Sets can be learned by equivalence queries in $O(n^2)$ rounds.

*Proof.* The Learner begins with the hypothesis $H$ of the complete graph. If the graph $G$ being taught is not complete, then the Teacher must respond with a positive counterexample $V'$. Since every singleton set is independent in any graph, $V'$ must contain at least two vertices. The Learner can deduce that there are no edges in $G$ between vertices in $V'$. The algorithm makes only equivalence queries. At each stage, an example of which is shown in Figure 3, the Learner presents a hypothesis graph $H$ that contains edges between all pairs of vertices except those pairs for which the Learner has deduced that no edge is present in $G$. It follows that the Teacher must respond with a positive counterexample, and this must allow for the non-presence of at least one new edge to be deduced. The algorithm will terminate in at

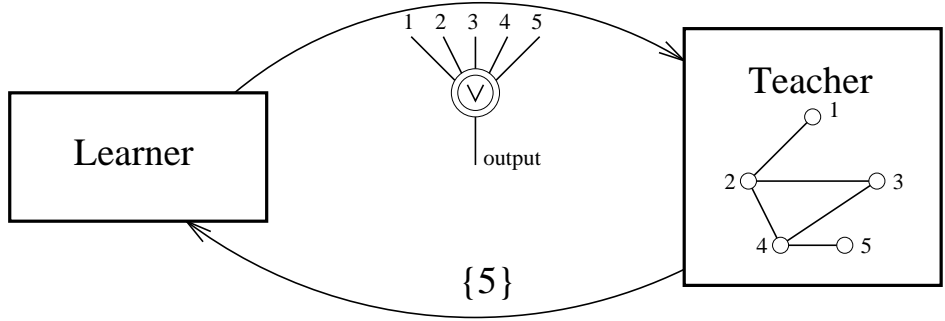

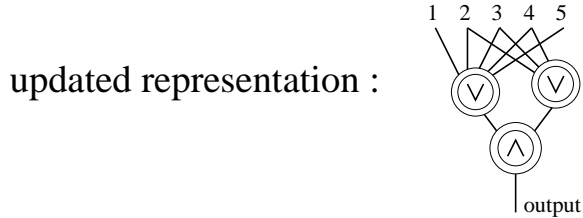

Figure 4: Update cycle for learning Dominating Sets by Extended Equivalence and Membership Queries.

most $\binom{n}{2}$ rounds. □

*Theorem 2.3* Dominating Sets can be learned by extended equivalence and membership queries in $O(n^2)$ rounds.

*Proof.* The Learner begins with a complete graph and makes an equivalence query. In the case where $G$ is not the complete graph, the Teacher must respond by providing a negative counterexample $V'$. The Learner now makes membership queries on supersets of $V'$ until a maximal negative counterexample is identified. The Learner can deduce from maximality that $V - V'$ consists precisely of the solid neighborhood $N[x]$ of some vertex $x$, $N[x] = \{u : u = x \text{ or } ux \in E(G)\}$. Note that the Learner will not know which of the vertices in $V - V'$ plays the role of $x$, however.

A hypothesis circuit is constructed (for an extended equivalence query) which "accepts" a set a vertices if and only if it has nonempty intersection with each solid neighborhood that has been identified. Inductively, if the hypothesis is not correct, then the Teacher must respond with a negative counterexample, which leads to the identification of another solid neighborhood. Figure 4 illustrates an example of the construction of an improved hypothesis. Since there are $n$ solid neighborhoods in the graph, and since a dominating set is equivalently a set having nonempty intersection with all solid neighborhoods, the Learner will produce a correct hypothesis circuit in $O(n^2)$ rounds. □

*Lemma 2.4* If $L$ is a learning algorithm for a property $\Phi$ of vertex sets that runs in $f(n)$ rounds and makes only equivalence queries, then $L$ also serves as a learning algorithm for the vertex sets having property $\Phi$ and cardinality $k$ that runs in at most $f(n)$ rounds.

*Proof.* The correctness of $L$ implies that in at most $f(n)$ rounds the Teacher will cease to respond with counterexamples for the predicate $\Phi$, regardless of the cardinality of the counterexamples that the teacher has been supplying. □

*Corollary.* $k$-Vertex Covers and $k$-Independent Sets can be exactly learned in $O(n^2)$ rounds by equivalence queries.

Our next theorem shows an improvement on the above for very small values of $k$. In considering the learning complexity of sets of vertices of a fixed cardinality $k$, it is important to keep in mind that this is precisely what is being taught, and nothing else. For example, we can make membership queries about sets of cardinality $k$, but the Teacher does not supply information about sets of other cardinalities. We will use the following structural lemma concerning minimal vertex covers in a graph.

*Lemma 2.5* If $G$ is a graph having a vertex cover of cardinality $m$, then $G$ has at most $2^m$ distinct minimal vertex covers.

*Proof.* Let $M$ denote the set of minimal vertex covers of $G$, and let $A, B \in M$, $A \neq B$. Note that $A - B \neq \emptyset \neq B - A$, since $A$ and $B$ are minimal covers.

Let $N(x) = \{y : xy \in E(G)\}$ and for $X \subseteq V(G)$, let $N(x, X) = N(x) \cap X$.

(1) For $x \in A - B$, $N(x) \subseteq B$, otherwise $B$ is not a vertex cover. It follows that $N(x, \overline{A}) = N(x, B - A) \subseteq B - A$. Similarly for $y \in B - A$.

(2) For $x \in A - B$, $N(x) \not\subseteq A$, otherwise $A$ is not a minimal vertex cover. It follows that $N(x, \overline{A}) = N(x, B - A) \neq \emptyset$. Similarly for $y \in B - A$.

(3) We have $\bigcup_{x \in A-B} N(x, \overline{A}) = B - A$ since for $y \in B - A$, there exist $x \in N(y, A)$, so that $y \in N(x, \overline{A})$. Similarly $\bigcup_{y \in B-A} N(y, \overline{B}) = A - B$.

Let $f : M \to \wp(A)$ by $f(B) = A \cap B$. Then $f$ is injective since $f(B_1) = f(B_2)$ implies $A \cap B_1 = A \cap B_2$. Consequently, $A - B_1 = A - B_2$ and from (3)

$$B_1 - A = \bigcup_{x \in A-B_1} N(x, \overline{A}) = \bigcup_{x \in A-B_2} N(x, \overline{A}) = B_2 - A$$

so that $B_1 = B_2$. Consequently $|M| \leq 2^m$. □

*Theorem 2.6* $k$-Vertex Covers can be learned by equivalence queries in $2^{k \cdot 2^k}$ rounds (independent of the order of the graph).

*Proof.* The learning algorithm can be represented as a tree of height at most $2^k$, with each node having at most $2^k$ children. Thus the number of nodes in the tree observes the bound in the statement of the theorem. A node in the tree represents the following sequence of activity:

(1) From a list of "presumed minimal vertex covers" received from the parent (initially empty at the root of the tree), construct a hypothesis graph $H$ containing all possible edges that have at least one endpoint in each vertex cover set on the received list (thus $H$ is the unique maximal graph consistent with the list). At the root node, in particular, $H$ is the complete graph.

(2) Make an equivalence query to the Teacher with the hypothesis graph constructed in (1).

• If the hypothesis graph $H$ is equivalent to $G$ with respect to $k$-element vertex covers, then of course we are done, and we will term this a *winning* node.

• If a negative counterexample is supplied by the Teacher, then the node is marked as "leaf" and has no descendants.

• If a positive counterexample $Q$ is supplied by the Teacher, then create one descendant node for each nonempty $C \subseteq Q$, corresponding to a "guess" that this subset of $Q$ is a newly discovered minimal vertex cover in $G$. We require that $C$ is nonempty because $\emptyset$ cannot be a minimal vertex cover of $G$. Since $Q$ has cardinality $k$, each node has no more than $2^k$ children. (Note that we cannot use membership queries to the Teacher to identify directly a subset of $Q$ that is a minimal vertex cover, because the Teacher supplies *only* information about sets of cardinality $k$.) Pass on to each descendant node the list received from the parent, augmented with the set of vertices in the corresponding subset of $Q$.

• If the list passed to any descendant contains more than $2^k$ distinctQJs2:s, then the descendant is identified as a *leaf* and is not processed further. The reason for this is that the presumption that the list contains only minimal vertex covers of $G$ is false, by Lemma 2.5.

Since only nodes that lead to positive counterexamples are expanded, and each positive counterexample leads to a larger list of presumed minimal vertex covers in each descendant, the tree has height at most $2^k$ by Lemma 2.5. This establishes that the tree is bounded in size as claimed.

The correctness of the algorithm follows from:

*Claim:* If at a node in the tree of the algorithm, the list received from the parent is a subset of the set of minimal vertex covers of $G$, then the hypothesis graph $H$ produced in step (1) is either equivalent to $G$, or admits only positive counterexamples.

To see this, suppose that $H$ and $G$ are inequivalent as witnessed by a negative counterexample $W$. Thus $W$ is a vertex cover for $H$ but not for $G$. The failure of $W$ to serve as a vertex cover in $G$ is witnessed by an edge $uv$ with $u \notin W$ and $v \notin W$. The edge $uv$ is absent in $H$, and by the construction of $H$, this can only be because there is a minimal vertex cover $M$ on the list that is incompatible with $uv$, i.e., $u \notin M$ and $v \notin M$. But our hypothesis is that $M$ is a (minimal) vertex cover in $G$, a contradiction.

Since the hypothesis of the above claim holds at the root trivially, and since the list passed to some descendant of an expanded node preserves the hypothesis (by representing a correct guess), we are done. □

## 3  Teaching Vertex Sets

We may consider also the complexity of *teaching* sets of vertices in the model of [GK91] and [JT92]. In this case, we assume that the Learner receives a set of positive examples and a set of negative examples from the Teacher. A *positive example* is a set of vertices $V'$ such that $\Phi(G, V')$ holds, and a *negative example* is a set of vertices $V'$ such that $\Phi(G, V')$ does not hold. A *teaching algorithm* $T$ for a predicate $\Phi$ is a strategy by which the Teacher can always lead the Learner to produce a graph $H$ that is $\Phi$-equivalent to the graph $G$ being taught, assuming that that the Learner is *competent*, that is, generates a hypothesis graph which is consistent with the examples that have previously been supplied by the Teacher. The teaching algorithm $T$ is said to run in $g(n)$ *rounds* if for any graph $G$ of order $n$ that is being taught, $\exists$ set of positive examples $P_G$ and set of negative examples $N_G$ where $\mid P_G \mid + \mid N_G \mid \leq g(n)$, such that any hypothesis $H$ where $\Phi(H, V_P)$ holds $\forall V_P \in P_G$ and $\Phi(H, V_N)$ holds $\forall V_N \in N_G$, is $\Phi$-equivalent to $G$.

*Theorem 3.1*  Vertex Covers can be taught in $O(n^2)$ rounds.

*Proof.*  For $G = (V, E)$, $u, v \in V$, an edge $uv \notin E$ iff $V - \{u, v\}$ is a vertex cover in $G$. The Teacher gives negative examples $V - \{u, v\}$ $\forall uv \in E$ and positive examples $V - \{u, v\}$ $\forall uv \notin E$. Any hypothesis $H = (V, E_H)$ consistent with these $\binom{n}{2}$ examples must have $E_H = E$, so $H = G$. □

*Theorem 3.2*  Independent Sets can be taught in $O(n^2)$ rounds.

*Proof.*  For $G = (V, E)$, $u, v \in V$, an edge $uv \notin E$ iff $\{u, v\}$ is an independent set in $G$. The Teacher gives negative examples $\{u, v\}$ $\forall uv \in E$ and positive examples $\{u, v\}$ $\forall uv \notin E$. Any hypothesis $H = (V, E_H)$ consistent with these $\binom{n}{2}$ examples must have $E_H = E$, so $H = G$. □

*Theorem 3.3* Teaching $k$-Independent Sets requires $\binom{n}{k}$ examples.

*Proof.* Consider the case $G = K_n$. There are no positive examples for $k \geq 2$. If $\mid N_G \mid < \binom{n}{k}$, $\exists$ negative example $V_N \notin N_G$. Let $H = (V, E_H)$ where $uv \in E_H$ iff $u \in V - V_N$ or $v \in V - V_N$. $\forall V' \in N_G$, $V'$ is not an independent set in $H$ since $\exists v \in V' - V_N$, so $\forall u \in V', u \neq v$, we have $uv \in E_H$. $H$ is therefore consistent with $N_G$. However, $\forall u, v \in V_N$, $uv \notin E_H$, so $V_N$ is an independent set in H but not in $K_n$. □

## 4 Summary and Open Problems

We believe that the results presented here on the learning complexity of some simple and familiar graph structures initiate an interesting chapter in the general study of learning algorithms and complexity, which is as yet a relatively new, though vigorously expanding research area.

The connections between the learnability of CNF formulas and the learnability of dominating sets is especially intriguing in view of the positive results obtained here concerning vertex covers and independent sets. The following questions would seem to be natural targets for further research.

(1) Is it possible to learn the $k$-element dominating sets in a graph by extended equivalence and membership queries in $o(n^k)$ rounds? (In other words, can we at least do better than the obvious "brute force" learning algorithm?)

(2) Can the sets of vertices that cover all of the cycles in a graph be learned in polynomial time? This would seem to be a property would be "harder" than vertex covers, but perhaps not as difficult as dominating sets.

(3) Is it possible to learn the $k$-element independent sets in a graph in a number of rounds $f(k)$ independent of the size of the graph?

(4) Can dominating sets be learned in polynomial time by (non-extended) equivalence and membership queries?

(5) What of lower bounds on the learning complexity of vertex sets?